\documentclass[12pt]{article}

\usepackage{latexsym}
\usepackage{setspace}
\usepackage{amsmath}
\usepackage{amssymb}
\usepackage{latexsym}
\usepackage{amsmath}
\usepackage{amssymb}
\usepackage{natbib}
\usepackage{subfigure}
\usepackage{color}
\setcitestyle{numbers}
\setcitestyle{square}
\usepackage{graphics}
\usepackage{graphicx}
\usepackage{multirow}
\makeatletter
\def\ps@pprintTitle{%
  \let\@oddhead\@empty
  \let\@evenhead\@empty
  \let\@oddfoot\@empty
  \let\@evenfoot\@oddfoot
}
\makeatother
\makeatletter \oddsidemargin  -.5cm \evensidemargin -.5cm \textwidth
17.5cm \topmargin -.65in \textheight 24cm

\newcommand{\bd}{\begin{definition}}
\newcommand{\ed}{\end{definition}}
\newcommand{\br}{\begin{remark}}
\newcommand{\er}{\end{remark}}
\newcommand{\bea}{\begin{eqnarray}}
\newcommand{\eea}{\end{eqnarray}}
\newcommand{\beann}{\begin{eqnarray*}}
\newcommand{\eeann}{\end{eqnarray*}}
\newtheorem{theorem}{Theorem}[section]

\newtheorem{lemma}[theorem]{Lemma}

\newtheorem{remark}{Remark}[section]
\newtheorem{example}{Example}[section]

\numberwithin{equation}{section}
\numberwithin{equation}{section}

\title{Inadmissibility of invariant estimator of function of scale parameter of several exponential distributions}
\author{\small $^a$Lakshmi Kanta Patra\footnote
	{\baselineskip=10pt
		~lkpatra@iitbhilai.ac.in, patralakshmi@gmail.com} ,$^a$Shrajal Bajpai, $^b$Neeraj Misra  \\
	\small $^a$Department of Mathematics, 
	\small Indian Institute of Technology Bhilai\\
    \small $^b$Department of Mathematics and Statistics, 
    \small Indian Institute of Technology Kanpur}
\spacing{1.1}
\begin{document}
\date{}
\maketitle
\begin{abstract}
In various applied areas such as reliability engineering, molecular biology, finance, etc., the measure of uncertainty of a probability distribution plays an important role. In the present work, we consider the estimation of a function of the scale parameter, namely entropy of many exponential distributions having unknown and unequal location parameters with a common scale parameter. For this estimation problem we have consider bowl-shaped location invariant loss functions. The inadmissibility of the minimum risk invariant estimator (MRIE) is proved by proposing a non-smooth improved estimator. Also, we have obtained a smooth estimator which improves upon the MRIE. As an application, we have obtained explicit expressions of improved estimators for two well known loss functions namely squared error loss and linex loss. Further we have shown that these estimators can be derived for other important censored sampling schemes. At first we have obtained the results for complete and i.i.d. sample. We have seen that the results can be applied for (i) record values, (ii) type-II censoring, and (iii) progressive Type-II censoring. Finally, a simulation study has been carried out to compare the risk performance of the proposed improved estimators.
\end{abstract}
\section{Introduction}
Similar to the hazard rate, entropy of a lifetime distribution is an important characteristic. It measures the uncertainty of a probability distribution. Shannon`s entropy is widely used in various areas of science and technology, such as ecology, hydrology, water resources, social studies, economics, biology, etc. In molecular sciences, estimation of the entropy of molecules plays an important role in understanding various chemical and biological processes \cite{nalewajski2002applications}. In economics, entropy estimation often allows the researchers to use data for the improvement of the assumptions on the parameters in econometric models, see \cite{golan1996maximum}. In reliability theory entropy is used
in to measure uncertainty \cite{kamavaram2002entropy}. If we want to estimate the uncertainty of a parallel or series system with several independent components, we need to predict uncertainty in individual components. A two-parameter experiential distribution is the most commonly used lifetime distribution in life testing experiments and reliability theory. In the case of the exponential distribution, entropy is a function of the scale parameter. The estimation of scale parameters and the function of scale parameters is a well-studied problem in statistical decision theory. The inadmissibility of the best affine invariant estimator of a normal variance was first established by \cite{stein1964}. This result motivates many statisticians to find improved estimators for scale parameters. One may refer to \cite{maatta1990developments} for a detailed review. The result of \cite{stein1964} was extended by \cite{brown1968inadmissibility} to prove the inadmissibility of the best equivariant estimator of powers of the scale parameter for a wide class of
location-scale densities and for a general invariant loss function. Two new techniques for obtaining improvements over equivariant estimators were developed by \cite{brewster1974alternative} for strictly bowl-shaped loss functions. In this paper, we have considered the estimation of a function of a common scale parameter, namely the entropy of several exponential distributions. Let $X$ be a random variable with
probability density function $f(x|\theta).$ Then, the Shannon's and R\'{e}nyi
entropy are given as
\begin{eqnarray}\label{eq0.1}
H(\theta)=-E(\ln f(x|\theta))
\end{eqnarray}
and
\begin{eqnarray}\label{eq1.1}
R_{\alpha}(\theta)=\frac{1}{1-\alpha}\ln\int_{-\infty}^{\infty}f^{\alpha}(x|\theta)dx=
\frac{1}{1-\alpha}\ln E\Big(f^{\alpha-1}(X|\theta)\Big)
\end{eqnarray}
respectively, where $\alpha\geq0$.\\

Entropy of a probability distribution is an important characteristic, like general moments, quantiles mean, median and standard deviation. Entropy gives us information about the uncertainty of probability distribution. Recent past, many authors have investigated the estimation of entropy of various probability models from a decision-theoretic point of view. Now we will describe some previous work in this direction. Estimation of entropy of a multivariate normal distribution has been considered by \cite{misra2005estimation}. They have shown the inadmissibility of usual estimators by deriving two improved estimators. Specially, they have obtained  Stein-type and Brewster-Zidek-type improved estimators. Finally, they have shown that the Brewster-Zidek-type improved estimator is generalized Bayes. Estimation of measure of uncertainty, that is, the entropy of several experiential distributions, was investigated by \cite{kayal2011estimating}. They have shown that the BAEE is inadmissible under squared error loss function. R\'{e}nyi entropy gives an important measure of uncertainty which is more flexible than Shannon entropy. Estimation  R\'{e}nyi entropy of $k$ exponential distributions with common location but different scale have been investigated by    \cite{kayal2015estimating}. The authors proposed a sufficient condition under which affine and scale equivariant estimators are inadmissible. \cite{patra2018estimatinga} discussed the problem of finding improved estimators of the entropy of an two-parameter exponential distribution with ordered location parameters. They have adopted the techniques of \cite{stein1964}, \cite{brewster1974improving}, and \cite{kubokawa1994unified} to find the improved estimators under a general location invariant loss function. \cite{patra2019estimating} proved that the usual estimator of the common hazard rate of several exponential distributions is inadmissible. They have obtained improved estimators which dominate the best affine equivariant estimator. Recently, \cite{patra2020improved} studied the problem of estimating the entropy of an exponential population based on a doubly censored sample. He proved the inadmissibility of the best affine equivariant estimator under a general bowl-shaped location invariant loss function. 

Let $\boldsymbol{X}_i=(X_{i1},\dots,X_{in})$ be a random sample taken from the population $\Pi_{i}$, $i=1,\dots,k$ $(k \ge2)$. We assume that the samples are taken independently. The population $\Pi
_{i}$ is assumed to have density
\begin{eqnarray}\label{model}
f_{i}(x;\theta_i,\sigma)=\left\{
        \begin{array}{ll}
           \frac{1}{\sigma}\exp\left(-\frac{x-\theta_i}{\sigma}\right),~~\mbox{ if }x>\theta_i,\theta_i \in \mathbb{R}, \sigma>0\\\\
           0,~~~~~~~~~~~~~~~~~~\mbox{ otherwise}
        \end{array}.
         \right.
\end{eqnarray}
For a population with probability density function (\ref{model}), the Shannon's entropy is $H(\sigma)=(1+\ln \sigma )$ and the Renyi entropy can be obtained as $R(\sigma)=\ln \sigma - \frac{\ln \alpha}{1-\alpha}$. So estimation of $H(\sigma)$ and $R(\sigma)$ is equivalent to estimation of $\theta=\ln \sigma$.\\

In this paper we have considered the estimation $\theta$ with respect to a location invariant loss function $L(T-\theta)$ where $L(t)$ statistics the following conditions.
\begin{itemize}
 \item[(i)] $L(t)$ is real valued absolute continuous and non monotone function. 
\item[(ii)] $L(t)$ is such that $L(t)$ is decreasing for $t<0$ and $L(t)$ is increasing for $t>0$ and $L(t)>0$ for all $t$.
 \end{itemize}
As a consequence of these conditions we get $L(t)$ is differentiable all most everywhere. Based on the $i$-th sample $(X_{i1},\dots,X_{in})$, a complete sufficient statistic for $(\theta_i,\sigma)$ is $(X_{i},Y_i)$ with $X_{i}=nX_{i(1)}$ and $Y_i=\sum_{j=1}^{n}(X_{ij}-X_{i(1)})$, where $X_{i(1)}=\min\{X_{i1},\dots,X_{in}\}$. We have $S=\sum_{i=1}^{k}Y_i$, then $V=\frac{S}{\sigma}$ follows a $Gamma(k(n-1),1)$ distribution.  Based on all the samples a complete sufficient statistic is $(\boldsymbol{X},S)$, where $\boldsymbol{X}=(X_1,\dots,X_k)$. So we have the pdf of $X_i$ and $V$ are

\begin{equation*}
h_i(x_i)=\frac{1}{\sigma}\exp \left\{-\frac{1}{\sigma}\left(x_i-n\theta_i\right)\right\},x_i>n\theta_i~~~~g(v)=\frac{e^{-v}v^{k(n-1)-1}}{\Gamma(k(n-1))},v>0
\end{equation*}
respectively.

Consider a group of transformation as $\mathcal{G}=\{g_{a,\boldsymbol{b}}, a>0,~\boldsymbol{b} \in \mathbb{R}^k\}$. The transformation is given as $g_{a,\boldsymbol{b}}(\boldsymbol{z}) =(a\boldsymbol{z}+\boldsymbol{b})$. Under this transformation complete sufficient statistics $(\boldsymbol{X},S)$ is invariant and we have 
\begin{eqnarray*}
	(\boldsymbol{\theta},\sigma)\rightarrow (a\boldsymbol{\theta}+\boldsymbol{b},a\sigma)~~~~\mbox{and}~~~~~(\boldsymbol{X},S)\rightarrow \left(a\boldsymbol{X}+\boldsymbol{b},aS\right).
\end{eqnarray*}
Consequently, we have 
$$\ln \sigma \rightarrow \ln \sigma +\ln a$$
Also we have the loss function $L(T-\theta)$ is invariant. Hence the form of an invariant estimator is obtained as
$$T_{c}=\ln S+c,$$
where $c$ is a real constant.

The following lemma gives the minimum risk invariant estimator (MRIE). We denote $\boldsymbol{\theta}=(\theta_1,\theta_2,\dots,\theta_k)$. Now onwards we use the notation $E_{\boldsymbol{\theta},1}$ for $E_{\boldsymbol{\theta},\sigma=1}$ and $E_{\boldsymbol{0},1}$ for $E_{\boldsymbol{\theta}=\boldsymbol{0},\sigma=1}$.
\begin{lemma}
The MRIE of $\theta$ with respect to a general location invariant
loss function $L(t)$ is
\begin{eqnarray}\label{baee}
T_{0}=\ln S+q_0,
\end{eqnarray}
where $q_0$ minimizes
\begin{eqnarray}\label{equ1}
E_{\boldsymbol{\theta},1}\left(L\left(\ln S+c\right)\right).
\end{eqnarray}
\end{lemma}
\noindent \textbf{Proof:} The proof is simple and hence omitted for the sake of brevity. \hfill\(\Box\)
\begin{example}
Let $L(t)=t^2$. Then, we have $c_{0}=-\psi(nk-k)$, where $\psi(.)$ denotes digamma function. So the MRIE of $\theta$ is $T_{01}=\ln S-\psi(kn-k)$.
\end{example}

\begin{example}
Let $L(t)=e^{at}-at-1$, $a\neq0.$ Then, using (\ref{equ1}), we
obtain $c_{0}=\frac{1}{a}\ln\left(\frac{\Gamma(nk-k)}{\Gamma(nk+a-k)}\right)$, where $a>k(1-n)$. In this case, the MRIE of $\theta$ is $T_{02}=\ln S+\frac{1}{a}\ln\left(\frac{\Gamma(nk-k)}{\Gamma(nk+a-k)}\right)$.
\end{example}

In the present work we aim to obtain estimators which improve upon the MRIE of $\ln \sigma$. We have studied this for four important sampling schemes: (i) complete and i.i.d. sample, (ii) record values, (iii) type-II censoring, and (iv) progressive Type-II censoring. Here we have adopted the techniques of \cite{stein1964} and \cite{brewster1974improving} for finding improved estimators. Several researchers have studied the problem of finding an improved estimator of a scale parameter in the presence of unknown location parameter using these techniques. For some nice applications of these techniques, we refer to \cite{misra2002smooth, misra2006some,petropoulos2010class,kayal2015estimating,tripathi2017estimating,petropoulos2017estimation,patra2019minimax,patra2018estimating,mani2020quantile}, and references therein.\\

The rest of the paper is organized as follows. In Section \ref{sec2},  we have proved the MRIE is inadmissible by deriving an improved estimator, which is not smooth based on the i.i.d sample. As an application, we have proposed improved estimators for the squared error loss and linex loss functions. A smooth estimator is derived in Section \ref{sec3}, which dominates the MRIE. For squared error loss and linex loss functions, we have derived the explicit expression of the smoothed improved estimators. A Bayes estimator has been given in Section \ref{sec4}. We have compared the risk performance of the improved estimators in Section \ref{sec5}. Section \ref{sec6} has considered special sampling schemes such as record values, type-II censoring, and progressive Type-II censoring. For each sampling scheme, improved estimators are obtained using the result for the i.i.d sampling scheme. A simulation study has been carried out for the record values. Finally, in Section \ref{sec7}, we have given concluding remarks.

\section{Inadmissibility of  MRIE}\label{sec2}
In the previous section, we obtained the MRIE of $\ln \sigma$. Now we will find an estimator which improves the MRIE under loss function $L(t)$. \cite{patra2019estimating} proved a similar result to find improved estimators of the common hazard rate of several exponential distributions. We have adopted the techniques similar to \cite{patra2019estimating}.  Also, this result extends a result of \cite{patra2018estimatinga} from one dimension to $k$ dimension. For proving the inadmissibility of MRIE, we consider the following class of scale invariant estimators. 
\begin{eqnarray}\label{se1}
	T_{\zeta}(\boldsymbol{X},S)=\ln S +\zeta(\boldsymbol{Z}),
\end{eqnarray}
where $\boldsymbol{Z}=(Z_1,\dots,Z_k)$, $Z_i=\frac{X_i}{S}$ and $\zeta$ is a measurable real valued function. The theorem below proves that the MRIE $T_0$ is inadmissible by proposing a non smooth dominating estimator.
\begin{theorem}\label{thmst}
	Suppose $p_0$ be the unique solution of the equation
	\begin{eqnarray}\label{equ2}
		E\left(L^{\prime}\left(\ln Y+p_0\right)\right)=0.
	\end{eqnarray}
	where $Y\sim Gamma(kn,1)$. Define an estimator
	\begin{eqnarray}\label{stein}
		T_{\zeta_0}(\boldsymbol{X},S)=\left\{
		\begin{array}{ll}
			\ln S+\min\left\{q_0,p_0+\ln\left(\sum_{i=1}^kZ_i+1\right)\right\},~~~Z_i>0, i=1,\dots,k\\\\
			\ln S+q_0, ~~~~~~~~~~~~~~~~~~~~~~~~~~~~\mbox{otherwise}
		\end{array},
		\right.
	\end{eqnarray}
	where $q_0$ is the unique solution of (\ref{equ1}). Then $T_{\zeta_0}(\boldsymbol{X},S)$ dominates the MRIE $T_0$ under a general location invariant loss function $L(t)$ under the condition
	\begin{eqnarray}\label{equ3}
		E\left[L^{\prime}\left(\ln V+p_0\right)\right]<0.
	\end{eqnarray}
\end{theorem}

\noindent\textbf{Proof:} It is easy to see that the risk of the estimator  $T_{\zeta}$ depends on the unknown parameter $(\boldsymbol{\theta},\sigma)$ only through $\theta_1/\sigma,\dots,\theta_k/\sigma$.  So we consider $\sigma=1$ without any loss of generality. We can easily write the the expected loss of $T_{\zeta}(\boldsymbol{X},S)$ as

\begin{eqnarray*}
R(\boldsymbol{\theta},T_{\zeta})=E_{\boldsymbol{\theta}}E_{\boldsymbol{\theta}}\left[L(\ln V+\zeta(W))\big|\boldsymbol{Z}=\boldsymbol{Z}\right].
\end{eqnarray*}

Let us denote the conditional risk as
\begin{eqnarray}\label{th1eq1}
R_1(\boldsymbol{\theta},q)=E_{\boldsymbol{\theta}}\left[L(\ln V+q)\big|\boldsymbol{Z}=\boldsymbol{z}\right].
\end{eqnarray}

 Suppose $\boldsymbol{z}=(z_i,\dots,z_k)$ be such that $z_i >0$, for all $i$, and there exist $j$ such that $\theta_j>0$ for $1 \le j \le k$. Define $\beta_i=n\theta_i$ and denote $\mu=\max_{1\le i\le k}\{\beta_i/z_{i}\}$. The conditional density of $V$ given $\boldsymbol{Z}=\boldsymbol{z}$ is
\begin{eqnarray}
f_{V}(v|\boldsymbol{Z}=\boldsymbol{z})\propto e^{-\left(\sum_{i=1}^{k}z_i+1\right)v}v^{kn-1},~~v>\mu.
\end{eqnarray}
Based on the assumption of the loss function, it is easy to see that $R_1(\boldsymbol{\theta},c)$ is strictly bowl shaped. Suppose $q_{\boldsymbol{\theta}}(\boldsymbol{Z})$ minimizes $R_1(\boldsymbol{\theta},c)$. Then $q_{\boldsymbol{\theta}}(\boldsymbol{Z})$ is the unique solution of 
\begin{eqnarray}\label{th1eq2}
E\left(L^{\prime}\left(\ln V+q_{\boldsymbol{\theta}}(\boldsymbol{Z})\right)\big|
\boldsymbol{Z}=\boldsymbol{z}\right)=0.
\end{eqnarray}
It can seen that 
\begin{eqnarray*}
E\left(L^{\prime}\left(\ln V+q_{\boldsymbol{\theta}}(\boldsymbol{Z})\right)\big| \boldsymbol{Z}=\boldsymbol{z}\right) >0,
\end{eqnarray*}
provided $e^{-q_{\boldsymbol{\theta}}(\boldsymbol{Z})} <\mu$,
which is a contradiction to (\ref{th1eq2}). From this we can conclude that
\begin{eqnarray}\label{th1eq3}
e^{-q_{\boldsymbol{\theta}}(\boldsymbol{Z})} > \mu.
\end{eqnarray}
Again we have $q_{\boldsymbol{0}}(\boldsymbol{Z})$ is the unique solution of
\begin{eqnarray}\label{th1eq4}
E_{\boldsymbol{\theta}=\boldsymbol{0}}\left(L^{\prime}\left(\ln V+q_{\boldsymbol{0}}(\boldsymbol{Z})\right)
\big|\boldsymbol{Z}=\boldsymbol{Z}\right)=0.
\end{eqnarray}
It can be seen that 
\begin{eqnarray}\label{th1eq5}
E_{\boldsymbol{\theta}=\boldsymbol{0}}\left(L^{\prime}\left(\ln V+q_{\boldsymbol{\theta}}(\boldsymbol{Z})\right)\big| \boldsymbol{Z}=\boldsymbol{Z}\right)<0.
\end{eqnarray}
for $q_{\boldsymbol{\theta}}(\boldsymbol{Z})>\mu$. 
Further we have $R_1(\boldsymbol{\theta},q_{\boldsymbol{\theta}}(\boldsymbol{Z}))$ is strictly bowl shaped. So from the equation (\ref{th1eq4}) and (\ref{th1eq5}), it is implied that 
\begin{eqnarray}\label{th1eq6}
q_{\boldsymbol{\theta}}(\boldsymbol{z})<q_{\boldsymbol{0}}(\boldsymbol{z}).
\end{eqnarray}

If $\theta_{i}\le 0$ for all $i=1,\dots,k$, then it is easy to see that
\begin{eqnarray*}
E_{\boldsymbol{\theta}}\left[L\left(\ln V+q\right)\big|\boldsymbol{Z}=
 \boldsymbol{z}\right]=E_{\boldsymbol{0}}\left[L\left(\ln V+q\right)\big|\boldsymbol{z}=
 \boldsymbol{z}\right]
\end{eqnarray*}
This enable us 
 $q_{\boldsymbol{\theta}}(\boldsymbol{z})=q_{\boldsymbol{0}}(\boldsymbol{z}).$
 By making a change of variable  $\mathfrak{u}=v(\sum_{i-1}^{k}z_i+1)$, we get from (\ref{th1eq4}) 
\begin{eqnarray}\label{th1eq7}
\int_{0}^{\infty}L^{\prime}\left(\ln \mathfrak{u}+q_{\boldsymbol{0}}(\boldsymbol{z})-\ln \left(\sum_{i=1}^{k}z_i+1\right)\right)
e^{-\mathfrak{u}}
\mathfrak{u}^{kn-1}d\mathfrak{u}=0.
\end{eqnarray}
From (\ref{th1eq7}) with (\ref{equ2}), we get $q_{\boldsymbol{0}}(\boldsymbol{z})=p_0+\ln\left(\sum_{i=1}^{k}z_i+1\right)$. Again (\ref{equ1}) and (\ref{equ3}) gives us
\begin{eqnarray}\label{th1eq8}
p_0<q_0.
\end{eqnarray}
Now we consider a function of the form
\begin{eqnarray}
\zeta_0(\boldsymbol{Z})=\left\{
        \begin{array}{ll}
        \min\left\{q_0,p_0+\ln\left(\sum_{i=1}^{k}z_i+1\right)\right\},~~~z_i>0,~i=1,\dots,k\\\\
        q_0,~~~~~~~~~~~~~~~~~~~~~~~~~~~~~~~\mbox{otherwise}
        \end{array}.
         \right.
\end{eqnarray}

From (\ref{th1eq6}) and (\ref{th1eq8}) we have $q_{\boldsymbol{\theta}}(\boldsymbol{Z})<\zeta_{0}(\boldsymbol{Z})<q_0$ on a set having probability grater than zero. We have  $R_1(\boldsymbol{\theta},q)$ is a strictly bowl shaped function in $q$,  $R_1(\boldsymbol{\theta},q)$ is a decreasing function of $q$ for $q_{\boldsymbol{\theta}}(\boldsymbol{z})<c$ and hence $R_1(\boldsymbol{\theta},\zeta_{0}(\boldsymbol{Z})) \le R_1(\boldsymbol{\theta},q_0)$. Since $\ln(\sum_{i=1}^{k}z_i+1) <q_0-p_0$ is satisfied on a set of positive probability for all $\boldsymbol{\theta}$ and $z_i>0$ for $i=1,\dots,k$. So obtained that Hence we get $$R(\boldsymbol{\theta},T_{\zeta_0}) \le R(\boldsymbol{\theta},T_{0}).$$
Hence the theorem is proved. \hfill\(\Box\)

\begin{example}
Let $L(t)=t^2$. Then, from (\ref{equ2}), we obtain $b_{0}=-\psi(kn)$.
Thus, the estimator
\begin{eqnarray*}
T_{01}^{*}=\left\{
        \begin{array}{ll}
        \ln S + \min\left\{-\psi(kn-k),\ln(1+\sum_{i=1}^kz_i)-\psi(kn)\right\},~~z_i>0,i=1,\dots,n\\
        \ln S - \psi(kn-k), ~~~~~~~~~~~~~~~~~~~~~~~~~~~~~~~~~~~\mbox{otherwise}
        \end{array}
         \right.
\end{eqnarray*}
dominates $T_{01}.$
\end{example}

\begin{example}
Let $L(t)=e^{at}-at-1,~a\neq0$. Then, from (\ref{equ2}), we obtain $b_{0}=\frac{1}{a}\ln\left(\frac{\Gamma(nk)}{\Gamma(a+nk)}\right)$. Thus, the estimator
\begin{eqnarray*}
T_{02}^{*}=\left\{
        \begin{array}{ll}
        \ln S +
        \min\left\{\frac{1}{a}\ln\left(\frac{\Gamma(nk-k)}{\Gamma(a+nk-k)}\right),
        \ln(1+\sum_{i=1}^kz_i)+\frac{1}{a}\ln\left(\frac{\Gamma(nk)}{\Gamma(a+nk)}\right)\right\},~~z_i>0\\
        \ln S+\frac{1}{a}\ln\left(\frac{\Gamma(nk-k)}{\Gamma(a+nk-k)}\right), ~~~~~
        ~~~~~~~~~~~~~~~~~~~~~~~~~~~~~~~~~~~~~~~\mbox{otherwise}
        \end{array}
         \right.
\end{eqnarray*}
dominates $T_{02}$ for $a>k(1-n)$.
\end{example}
\section{Brewster-Zidek type improved estimator}\label{sec3}
In the last section, we proposed a non-smooth estimator. In this section, we prove inadmissibility of MRIE  $T_0$ by proposing an estimators which is smooth for $\{\boldsymbol{Z}=(Z_1,\dots,Z_k):Z_i \in (0,\infty),i=1,2,\dots,k\}$. \cite{patra2019estimating} studied the problem of finding a smooth, improved estimator for common hazard rate. Here we consider finding smooth estimator of the entropy of several exponential distributions. For this purpose consider an estimator of the form
\begin{eqnarray}\label{bz1}
T_{d}(\boldsymbol{X},S)=\left\{
        \begin{array}{ll}
        \ln S+d,~~~\boldsymbol{Z}\in \mathcal{B}_{\boldsymbol{r}}\\\\
        \ln S+q_0 ~~~~~\mbox{otherwise}
        \end{array},
         \right.
\end{eqnarray}
where $\mathcal{B}_{\boldsymbol{r}}= (0,r_1] \times (0,r_2] \times \dots (0,r_k]$ with $r_i>0$ for all $i=1,\dots, k$. With out loss of generality again, we take $\sigma=1$. Define $\boldsymbol{\beta}=(\beta_1,\dots,\beta_k)$, where $\beta_i=n\theta_i$ and $\eta=\max_{1\le i \le k}\{\beta_i/r_i\}$. The results in this section extends a results of \cite{patra2018estimatinga} form single exponential distribution to several exponential distribution. 

To propose an improve estimators we will analyze the conditional risk function
\begin{eqnarray*}
\mathcal{F}(d,\boldsymbol{\beta})=E_{\boldsymbol{\beta}}\left[L\left(\ln S+d\right)\big|\boldsymbol{Z}
\in\mathcal{B}_{\boldsymbol{r}}\right]
\propto\int_{0}^{\infty}L(\ln v+d)f_{\boldsymbol{\beta}}(v,\boldsymbol{r})dv,
\end{eqnarray*}
where
\begin{eqnarray*}
f_{\boldsymbol{\beta}}(v,\boldsymbol{r})\propto v^{kn-k-1}e^{-v}
\prod_{i=1}^{k}\left(e^{\beta_iI_{(\beta_i>0)}}-e^{-vr_i}\right)I_{(v>\eta)}.
\end{eqnarray*}
The following lemma we will study the properties of conditional risk function. This is useful to prove inadmissibility of MRIE
\begin{lemma}
\begin{itemize}
\item[(1)] For every $r_i>0$, $i=1,\dots,k$ the conditional risk $\mathcal{F}(d,\boldsymbol{\beta})$ is a strictly bowl shaped function in $d$. 
\item[(2)] Suppose $d(\boldsymbol{r}$ $ \boldsymbol{\beta})$  be the minimizer of $\mathcal{F}(d,\boldsymbol{\beta})$ and $d(\boldsymbol{r},\boldsymbol{0})$ be the minimizers of  $\mathcal{F}(d,\boldsymbol{0})$. Then for all $\boldsymbol{\beta}$ we have $d(\boldsymbol{r},\boldsymbol{\beta}) \le d(\boldsymbol{r},\boldsymbol{0}).$ 
\item[(3)] The function $d(\boldsymbol{r},\boldsymbol{0})$ is non decreasing in $r_i$ for $i=1,\dots,k$.
\end{itemize}
\end{lemma}
\noindent\textbf{Proof:}~~\textbf{(1)} To prove $\mathcal{F}(d,\boldsymbol{\beta})$ is a strictly bowl shaped by Lemma (2.1) of \cite{brewster1974improving} we have to prove that $\frac{\Lambda_{\boldsymbol{\beta}}(y-d_2,\boldsymbol{r})}
{\Lambda_{\boldsymbol{\beta}}(y-d_1,\boldsymbol{r})}$ is increasing in $y$ for given $0<d_1<d_2$, where
\begin{eqnarray*}
\Lambda_{\boldsymbol{\beta}}(y,\boldsymbol{r})\propto (e^y)^{kn-k}e^{-e^y}
\prod_{i=1}^{k}\left(e^{\beta_iI_{(\beta_i>0)}}-e^{-r_ie^y}\right)I_{(e^y>\eta/r)}.
\end{eqnarray*}
Now we have
\begin{eqnarray*}
\frac{\Lambda_{\boldsymbol{\beta}}(y-d_2,\boldsymbol{r})}
{\Lambda_{\boldsymbol{\beta}}(y-d_1,\boldsymbol{r})}=
\frac{(e^{y-d_2})^{kn-k}e^{-e^{y-d_2}}
\prod_{i=1}^{k}\left(e^{\beta_iI_{(\beta_i>0)}}-e^{-r_ie^{y-d_2}}\right)I_{(e^{y-d_2}>\eta/r)}}
{(e^{y-d_1})^{kn-k}e^{-e^{y-d_1}}
\prod_{i=1}^{k}\left(e^{\beta_iI_{(\beta_i>0)}}-e^{-r_ie^{y-d_1}}\right)I_{(e^{y-d_1}>\eta/r)}}.
\end{eqnarray*}
Using Lemma 6.1 of \cite{patra2018estimatinga} it can be easily proved that $\frac{\Lambda_{\boldsymbol{\beta}}(y-d_2,\boldsymbol{r})}
{\Lambda_{\boldsymbol{\beta}}(y-d_1,\boldsymbol{r})}$ is increasing in $y$ for $0<d_1<d_2$.\\

\noindent\textbf{(2)} By part (1) we have $\mathcal{F}(d,\boldsymbol{\beta})$ is a strictly bowl shaped function in $d$. Hence we get unique minimizer and we have
\begin{eqnarray}\label{sme1}
\mathcal{F}^{\prime}(d(\boldsymbol{r},\boldsymbol{\beta}),\boldsymbol{\beta})=0,~~\mbox{for all~} \boldsymbol{\beta}.
\end{eqnarray}
Suppose $d(\boldsymbol{r},\boldsymbol{\beta}) > d(\boldsymbol{r},\boldsymbol{0})$, then we have
\begin{eqnarray}
\mathcal{F}^{\prime}(d(\boldsymbol{r},\boldsymbol{\beta}),\boldsymbol{\beta})>
\mathcal{F}^{\prime}(d(\boldsymbol{r},\boldsymbol{0}),\boldsymbol{\beta}))>
\mathcal{F}^{\prime}(d(\boldsymbol{r},\boldsymbol{0}),\boldsymbol{0}))=0
\end{eqnarray}
which is a contradiction to (\ref{sme1}). This proves that $d(\boldsymbol{r},\boldsymbol{\beta}) \le d(\boldsymbol{r},\boldsymbol{0})$  and the last inequality in the above expression follows form the fact that $\frac{\Lambda_{\boldsymbol{\beta}}(y,\boldsymbol{r})}
{\Lambda_{\boldsymbol{0}}(y,\boldsymbol{r})}$ is increasing in $y$.\\

\noindent\textbf{(3)} Let $r_i=t$ and denote $d(\boldsymbol{r_t},\boldsymbol{0})=d((r_1,\dots,r_{i-1},t,r_{i+1},\dots,r_k),\boldsymbol{0})$. Now for  $0<t<t_1$ by Lemma 6.2 of \cite{patra2018estimatinga} it can be easily seen that $\frac{\Lambda_{\boldsymbol{0}}(y,\boldsymbol{r_t})}
{\Lambda_{\boldsymbol{0}}(y,\boldsymbol{r_{t_1}})}$ is nondecreasing in $y$ which implies that $d(\boldsymbol{r_t},\boldsymbol{0})$ is non decreasing in $t$. This proves that $d(\boldsymbol{r},\boldsymbol{0})$ is non decreasing in $r_i$ for $i=1,\dots,k$.

\hfill $\square$

As a consequence of of the above lemma we have the following dominance result. 
\begin{theorem}
The estimator $T_{d(\boldsymbol{r},\boldsymbol{0})}=\zeta_{\boldsymbol{r}}(\boldsymbol{Z})+\ln S$ dominates $T_0$ under a general location invariant loss function $L(t)$, where
\begin{eqnarray}
\zeta_{\boldsymbol{r}}(\boldsymbol{Z})=\left\{
        \begin{array}{ll}
        d(\boldsymbol{r},\boldsymbol{0}),~~~\boldsymbol{Z} \in \mathcal{B}_{\boldsymbol{r}}\\\\
        q_0,  ~~~~~~\mbox{otherwise}.
        \end{array}
         \right.
\end{eqnarray}

\end{theorem}
%
Let consider $\boldsymbol{r}'=(r'_1,\dots,r'_k)$ be such that $r'_i<r_i$ for $i=1,\dots,k$ and denote $\mathcal{B}_{\boldsymbol{r}',\boldsymbol{r}}=(r'_1,r_1]\times \dots (r'_k,r_k]$. Define an estimator as below. 
\begin{eqnarray*}\label{bz2}
T_{d,\boldsymbol{r},\boldsymbol{r}'}(\boldsymbol{X},S)=\left\{
        \begin{array}{ll}
        \ln S+d,~~~~~~~~~~\boldsymbol{Z}\in \mathcal{B}_{\boldsymbol{r}'}\\\\
         \ln S+d(\boldsymbol{r},\boldsymbol{0}),~~~\boldsymbol{Z}\in \mathcal{B}_{\boldsymbol{r}',\boldsymbol{r}}\\\\
        \ln S+q_0, ~~~~~~~~~~\mbox{otherwise}.
        \end{array},
         \right.
\end{eqnarray*}
Proceeding as above it we can easily prove that the estimator $T_{c(\boldsymbol{r}',\boldsymbol{0}),\boldsymbol{r},\boldsymbol{r}'}=\ln S+\zeta_{\boldsymbol{r},
\boldsymbol{r}'}(\boldsymbol{Z})$ dominates
$T_{c(\boldsymbol{r},\boldsymbol{0})}$ and hence the MRIE $T_{0}$, where
\begin{eqnarray*}
\zeta_{\boldsymbol{r},\boldsymbol{r}'}(\boldsymbol{Z})=\left\{
        \begin{array}{ll}
        d(\boldsymbol{r}',\boldsymbol{0}),~~~~ \boldsymbol{Z} \in \mathcal{B}_{\boldsymbol{r}'}\\\\
        d(\boldsymbol{r},\boldsymbol{0}),~~~\boldsymbol{Z}\in\mathcal{B}_{\boldsymbol{r}',
        \boldsymbol{r}}\\\\
        q_0,  ~~~~~~~~~\mbox{otherwise}
        \end{array}.
         \right.
\end{eqnarray*}

Similar to \cite{patra2019estimating} now we select a partition of $\prod_{i=1}^k[0,\infty)$ as $0=r^i_{j,0}<r^i_{j,1}<r^i_{j,2}<\dots<r^i_{j,m_j-1}<r^i_{j,m_j}<\infty$ for $i=1,2,\dots,k$.
\begin{eqnarray*}
\left\{
        \begin{array}{ll}
       0=r^1_{j,0}<r^1_{j,1}<r^1_{j,2}<\dots<r^1_{j,m_j-1}<r^1_{j,m_j}<\infty\\\\
       0=r^2_{j,0}<r^2_{j_1}<r^2_{j,2}<\dots<r^2_{j,m_j-1}<r^2_{j,m_j}<\infty\\\\
       0=r^3_{j,0}<r^3_{j,1}<r^3_{j,2}<\dots<r^3_{j,m_j-1}<r^3_{j,m_j}<\infty\\
       .\\
       .\\
       .\\
       0=r^k_{j,0}<r^k_{j,1}<r^k_{j,2}<\dots<r^k_{j,m_j-1}<r^k_{j,m_j}<\infty
        \end{array}
         \right.
\end{eqnarray*}
Let $\boldsymbol{r}_{j,l}=(r^1_{j,l},\dots,r^k_{j,l})$, $l=0,\dots,m_j$ and define
\begin{eqnarray*}
\zeta_j(\boldsymbol{Z})=\left\{
        \begin{array}{ll}
        d\left(\boldsymbol{r}_{j,1},\boldsymbol{0}\right),~~~\mbox{ if }~~r^i_{j,0}<z_i\le r^i_{j,1}, i=1,\dots,k\\\\
        d\left(\boldsymbol{r}_{j,2},\boldsymbol{0}\right),~~~\mbox{ if }~~r^i_{j,1}<z_i\le r^i_{j,2}, i=1,\dots,k\\\\
        d\left(\boldsymbol{r}_{j,3},\boldsymbol{0}\right),~~~\mbox{ if }~~r^i_{j,2}<z_i\le r^i_{j,3}, i=1,\dots,k\\\\
        .\\
        .\\
        d\left(\boldsymbol{r}_{j,m_j},\boldsymbol{0}\right),~~\mbox{ if }~~r^i_{j,m_j-1}<z_i\le r^i_{j,m_j}, i=1,\dots,k\\\\
        q_0,~~~~~~~~~~~~~~~~~~~~~~~~\mbox{ otherwise}
        \end{array}
         \right.
\end{eqnarray*}
Assume that
\begin{eqnarray*}
\max_{1\le \kappa\le m_j}|r^{\nu}_{j,\kappa}-r^{\nu}_{j,\kappa-1}|\rightarrow 0 \mbox{ as } j\rightarrow\infty \mbox{~for~} \nu=1,\dots,k .
\end{eqnarray*}
Then we have $\zeta_j(\boldsymbol{Z}) \rightarrow \zeta_{*}(\boldsymbol{Z}) \mbox{ pointwise as } j \rightarrow \infty.$
Since for $j=1,2,\dots$, the estimator $T_{j}(\boldsymbol{X},S)=\ln S+\zeta_j(\boldsymbol{Z})$ has smaller risk than that of $T_0$. Then by applying Fatou's lemma we have the following dominance result.
\begin{theorem}\label{bzthm}
Define a function of the form 
\begin{eqnarray}\label{bz3}
\zeta_{*}(\boldsymbol{Z})=\left\{
        \begin{array}{ll}
        d(\boldsymbol{Z},\boldsymbol{0}),~~~~z_1>0,\dots,z_k>0\\\\
        q_0,  ~~~~~~~~~\mbox{otherwise}
        \end{array}.
         \right.
\end{eqnarray}
Then he estimator $T_{BZ}(\boldsymbol{X},S)=\ln S+\zeta_{*}(\boldsymbol{Z})$ dominates $T_0$ with respect to a general location invariant loss function $L(t)$.
\end{theorem}
Now we consider two special loss function and derive the smooth improved estimators which have uniformly smaller risk than $T_0$. 
\begin{example}\label{examp1}
	Consider the squared error loss function $L(t)=t^2$. Then from Theorem \ref{bzthm} we get the Brewster-Zidek type estimator as
	\begin{eqnarray}
		T_{BZ1}(\boldsymbol{X},S)=\left\{
		\begin{array}{ll}
			\ln S+d(\boldsymbol{Z},\boldsymbol{0}),~~~~z_1>0,\dots,z_k>0\\\\
			\ln S-\psi(nk-k),  ~~~~~~~~~\mbox{otherwise}
		\end{array}.
		\right.
	\end{eqnarray}
where 
\begin{eqnarray*}
	d(\boldsymbol{z},\boldsymbol{0})=-\frac{\displaystyle \int_{0}^{\infty}\ln v e^{-v}v^{kn-k-1}\prod_{i=1}^k(1-e^{-vz_i})dv}{\displaystyle\int_{0}^{\infty}e^{-v}v^{kn-k-1}\prod_{i=1}^k(1-e^{-vz_i})dv}
\end{eqnarray*}
For $k=2$
\begin{eqnarray*}
	d(\boldsymbol{z},\boldsymbol{0})=-\frac{\psi(2n-2)\left[1-(z_1+1)^{2-2n}-(z_2+1)^{2-2n}+(z_1+z_2+1)^{2-2n}\right]+D}{1-(z_1+1)^{2-2n}-(z_2+1)^{2-2n}+(z_1+z_2+1)^{2-2n}},
\end{eqnarray*}
where 
\begin{eqnarray*}
	D=\frac{\ln (1+z_1)}{(1+z_1)^{2n-2}}+\frac{\ln (1+z_2)}{(1+z_2)^{2n-2}}-\frac{\ln (1+z_1+z_2)}{(1+z_1+z_2)^{2n-2}}
\end{eqnarray*}
\end{example}
\begin{example}\label{examp2}
	Consider the linex loss function $L(t)=e^{at}-at-1, a\ne 0$. Then from Theorem \ref{bzthm} we get the Brewster-Zidek type estimator as
	\begin{eqnarray}
		T_{BZ2}(\boldsymbol{X},S)=\left\{
		\begin{array}{ll}
			\ln S+d(\boldsymbol{Z},\boldsymbol{0}),~~~~Z_1>0,\dots,Z_k>0\\\\
			\ln S+\frac{1}{a}\ln \left(\frac{\Gamma(nk-k)}{\Gamma(nk+a-k)}\right) ,  ~~~~~~~~~\mbox{otherwise}
		\end{array}.
		\right.
	\end{eqnarray}
	where $d(\boldsymbol{z},\boldsymbol{0})=-\frac{1}{a}\ln H(\boldsymbol{z})$ with 
	\begin{eqnarray*}
		H(\boldsymbol{z})=\displaystyle\int_{0}^{\infty}e^{-v}v^{kn+a-k-1}\prod_{i=1}^k(1-e^{-vz_i})dv
	\end{eqnarray*}
For $k=2$
$$H(\boldsymbol{z})=\Gamma(a+2n-2)\left(1-(z_1+1)^{2-2n-a}-(z_2+1)^{2-2n-a}+(z_1+z_2+1)^{2-2n-a}\right)$$
\end{example}
\section{Bayes estimator}\label{sec4}
In this section, we have consider the Bayes estimation of $\theta=\ln \sigma$. For this purpose we consider the prior distribution as
\begin{eqnarray}
\Pi(\boldsymbol{\theta},\sigma)=\left(\prod_{i=1}^k\frac{1}{\sigma}\exp\left(\frac{\theta_i-\mu_0}{\sigma}
\right)I(\theta_i<\mu_0)\right)\frac{\sigma_0^{\nu+1}}{\Gamma(\nu)\sigma^{\nu+1}}
\exp\left(-\frac{\sigma_0}{\sigma}\right).
\end{eqnarray}
So we have the posterior distribution of $(\boldsymbol{\theta},\sigma)$ is obtained as
\begin{eqnarray*}
\Pi(\boldsymbol{\theta},\sigma|\boldsymbol{X})=\left(\prod_{i=1}^k\frac{1+n_i}{\sigma}\exp
\left(\frac{n_i+1}{\sigma}(\theta_i-\min\{x_{i(1)},\mu_0\})
\right)I(\theta_i<\min\{x_{i(1)},\mu_0\})\right)\\
\frac{1}{\sigma^{\nu+1+nk}}
\exp\left(-\frac{1}{\sigma}(k\mu_0+\sigma_0+\sum_i\sum_jx_{ij})\right).
\end{eqnarray*}
Now we have for given $\sigma$; $\theta_1,\theta_2,\dots,\theta_k$ are independent with $\frac{n_i+1}{\sigma}(\theta_i-\min\{x_{i(1)},\mu_0\}) \sim Exp(1)$. Also we get $\frac{1}{\sigma}(k\mu_0+\sigma_0+\sum_i\sum_jx_{ij}) \sim \Gamma(nk+\nu)$. \\

So the the Bayes estimator with respect to squared error loss function is obtained as
\begin{eqnarray*}
T_{\Pi}=E(\ln \sigma|\boldsymbol{X})=\ln\left(k\mu_0+\sigma_0+\sum_i\sum_jx_{ij}\right)-\zeta(nk+\nu).
\end{eqnarray*}
\section{Simulation study}\label{sec5}
In the above, we have discussed the inadmissibility of MIRE of $\ln \sigma$. To prove the inadmissibility, we have obtained two improved estimators. One is Stein-type non smooth estimators, and the other one is Brewster-Zidek-type smooth improved estimator of entropy of several exponential distributions. In this section, we will study the risk performance of the proposed improved estimators with respect to the squared error loss function by simulation. For the purpose of simulation, we have generated 20,000 random samples from two exponential distributions with location parameters $\theta_1$, $\theta_2$, and scale parameter $1$ of sizes $n=4,6,8$. In Table \ref{table1}, we have tabulated percentage risk improvement (PRI) with respect to the MRIE $T_0$ of the improved estimators for $n=4$ and $n=6$ and $n=8$, respectively. For the simulation study, we have considered different values of $\theta_1$ and $\theta_2$. The PRI of an estimator $T$  with respect to $T_0$ is defined as
\begin{eqnarray*}
	PRI(T)=\frac{Risk(T_0)-Risk(T)}{Risk(T_0)} \times 100
\end{eqnarray*}
From the simulation, we have the following observations.
\begin{itemize}
	\item [(i)] The PRI of $T_{01}^{*}$ decreases as the value of $\theta_1$ and $\theta_2$ increases. The risk performance of $T_{01}^{*}$ is better than  $T_{BZ1}$ when the values of $\theta_1$ and $\theta_2$ near zero. 
	\item[(ii)] The PRI of $T_{BZ1}$ increases and then decreases as the values of $\theta_1$ and $\theta_2$ increases. The risk performance of  $T_{BZ1}$ better than  $T_{01}^{*}$ for larger values of $\theta_1$ and $\theta_2$. 
	\item [(iii)] We have seen that as $n$ increases the PRI of $T_{01}^{*}$ and $T_{BZ1}$ decrease. For large values of $n$, the performance of $T_{01}^{*}$ and the MIRE are the same. For large values of risk performance of $T_{BZ1}$ better than  $T_{01}^{*}$.
	\item [(iv)] Form the simulation we say that overall $T_{BZ1}$ perform better than the other estimators. 
\end{itemize}
We have similar type of observations for linex loss function. 
\begin{table}[h!]
	\centering
	\caption{Percentage risk improvement with respect to squared error loss function}
	\vspace{0.2cm}
	\label{table1}
	\begin{tabular}{|cc|cc|cc|cc|}
		\hline
		\multicolumn{2}{|c|}{$n$}                                & \multicolumn{2}{c|}{$4$}                     & \multicolumn{2}{c|}{$6$}                     & \multicolumn{2}{c|}{$8$}                    \\ \hline
		\multicolumn{1}{|c|}{$\theta_2$}                & $\theta_1$ & \multicolumn{1}{c|}{$T_{01}^{*}$}          & $T_{BZ1}$      & \multicolumn{1}{c|}{$T_{01}^{*}$}          & $T_{BZ1}$      & \multicolumn{1}{c|}{$T_{01}^{*}$}         & $T_{BZ1}$      \\ \hline
		\multicolumn{1}{|c|}{\multirow{6}{*}{$0.1$}} & 0.1     & \multicolumn{1}{c|}{10.47579}     & 0.665244 & \multicolumn{1}{c|}{5.23372}      & 2.720082 & \multicolumn{1}{c|}{2.196369}    & 3.563242 \\ \cline{2-8} 
		\multicolumn{1}{|c|}{}                       & 0.2     & \multicolumn{1}{c|}{8.177379}     & 3.837873 & \multicolumn{1}{c|}{2.352601}     & 4.995069 & \multicolumn{1}{c|}{0.370062}    & 5.044513 \\ \cline{2-8} 
		\multicolumn{1}{|c|}{}                       & 0.5     & \multicolumn{1}{c|}{2.071238}     & 6.912677 & \multicolumn{1}{c|}{0.02876744}   & 5.791525 & \multicolumn{1}{c|}{0}           & 4.552717 \\ \cline{2-8} 
		\multicolumn{1}{|c|}{}                       & 0.6     & \multicolumn{1}{c|}{1.128674}     & 7.042573 & \multicolumn{1}{c|}{0.0055315}  & 5.500053 & \multicolumn{1}{c|}{0}           & 4.190413 \\ \cline{2-8} 
		\multicolumn{1}{|c|}{}                       & 0.7     & \multicolumn{1}{c|}{0.564448}     & 7.004393 & \multicolumn{1}{c|}{0.0004631} & 5.185024 & \multicolumn{1}{c|}{0}           & 3.894156 \\ \cline{2-8} 
		\multicolumn{1}{|c|}{}                       & 0.8     & \multicolumn{1}{c|}{0.2685096}    & 6.872046 & \multicolumn{1}{c|}{0}            & 4.888716 & \multicolumn{1}{c|}{0}           & 3.665736 \\ \hline
		\multicolumn{1}{|c|}{\multirow{6}{*}{0.2}}   & 0.1     & \multicolumn{1}{c|}{8.177379}     & 3.888181 & \multicolumn{1}{c|}{2.352601}     & 5.0272   & \multicolumn{1}{c|}{0.370062}    & 5.054616 \\ \cline{2-8} 
		\multicolumn{1}{|c|}{}                       & 0.2     & \multicolumn{1}{c|}{5.611596}     & 6.894204 & \multicolumn{1}{c|}{0.7650419}    & 7.086941 & \multicolumn{1}{c|}{0.02789947}  & 6.307205 \\ \cline{2-8} 
		\multicolumn{1}{|c|}{}                       & 0.5     & \multicolumn{1}{c|}{1.128674}     & 9.681164 & \multicolumn{1}{c|}{0.0055315}  & 7.550914 & \multicolumn{1}{c|}{0}           & 5.507273 \\ \cline{2-8} 
		\multicolumn{1}{|c|}{}                       & 0.6     & \multicolumn{1}{c|}{0.564448}     & 9.752206 & \multicolumn{1}{c|}{0.0004630} & 7.200807 & \multicolumn{1}{c|}{0}           & 5.099963 \\ \cline{2-8} 
		\multicolumn{1}{|c|}{}                       & 0.7     & \multicolumn{1}{c|}{0.2685096}    & 9.665549 & \multicolumn{1}{c|}{0}            & 6.840881 & \multicolumn{1}{c|}{0}           & 4.772228 \\ \cline{2-8} 
		\multicolumn{1}{|c|}{}                       & 0.8     & \multicolumn{1}{c|}{0.1268464}    & 9.492634 & \multicolumn{1}{c|}{0}            & 6.509755 & \multicolumn{1}{c|}{0}           & 4.521561 \\ \hline
		\multicolumn{1}{|c|}{\multirow{6}{*}{0.4}}   & 0.1     & \multicolumn{1}{c|}{3.523156}     & 6.578196 & \multicolumn{1}{c|}{0.1702922}    & 6.007922 & \multicolumn{1}{c|}{0.001014834} & 4.944141 \\ \cline{2-8} 
		\multicolumn{1}{|c|}{}                       & 0.2     & \multicolumn{1}{c|}{2.071238}     & 9.369283 & \multicolumn{1}{c|}{0.02876744}   & 7.813017 & \multicolumn{1}{c|}{0}           & 5.953738 \\ \cline{2-8} 
		\multicolumn{1}{|c|}{}                       & 0.5     & \multicolumn{1}{c|}{0.2685096}    & 11.74199 & \multicolumn{1}{c|}{0}            & 7.835906 & \multicolumn{1}{c|}{0}           & 4.784317 \\ \cline{2-8} 
		\multicolumn{1}{|c|}{}                       & 0.6     & \multicolumn{1}{c|}{0.1268464}    & 11.72101 & \multicolumn{1}{c|}{0}            & 7.400792 & \multicolumn{1}{c|}{0}           & 4.317633 \\ \cline{2-8} 
		\multicolumn{1}{|c|}{}                       & 0.7     & \multicolumn{1}{c|}{0.05596951}   & 11.5567  & \multicolumn{1}{c|}{0}            & 6.974144 & \multicolumn{1}{c|}{0}           & 3.947277 \\ \cline{2-8} 
		\multicolumn{1}{|c|}{}                       & 0.8     & \multicolumn{1}{c|}{0.02279446}   & 11.31756 & \multicolumn{1}{c|}{0}            & 6.59024  & \multicolumn{1}{c|}{0}           & 3.665847 \\ \hline
		\multicolumn{1}{|c|}{\multirow{6}{*}{0.7}}   & 0.1     & \multicolumn{1}{c|}{0.564448}     & 7.089663 & \multicolumn{1}{c|}{0.0004630} & 5.211967 & \multicolumn{1}{c|}{0}           & 4.476075 \\ \cline{2-8} 
		\multicolumn{1}{|c|}{}                       & 0.2     & \multicolumn{1}{c|}{0.2685096}    & 9.700335 & \multicolumn{1}{c|}{0}            & 6.835372 & \multicolumn{1}{c|}{0}           & 5.686527 \\ \cline{2-8} 
		\multicolumn{1}{|c|}{}                       & 0.5     & \multicolumn{1}{c|}{0.02279446}   & 11.68276 & \multicolumn{1}{c|}{0}            & 6.5038   & \multicolumn{1}{c|}{0}           & 4.685471 \\ \cline{2-8} 
		\multicolumn{1}{|c|}{}                       & 0.6     & \multicolumn{1}{c|}{0.009141}  & 11.5671  & \multicolumn{1}{c|}{0}            & 5.993729 & \multicolumn{1}{c|}{0}           & 4.14633  \\ \cline{2-8} 
		\multicolumn{1}{|c|}{}                       & 0.7     & \multicolumn{1}{c|}{0.00295038}   & 11.32069 & \multicolumn{1}{c|}{0}            & 5.506559 & \multicolumn{1}{c|}{0}           & 3.682823 \\ \cline{2-8} 
		\multicolumn{1}{|c|}{}                       & 0.8     & \multicolumn{1}{c|}{0.0008821} & 11.00996 & \multicolumn{1}{c|}{0}            & 5.073649 & \multicolumn{1}{c|}{0}           & 3.304135 \\ \hline
	\end{tabular}
\end{table}

\section{Special sampling schemes}\label{sec6}
In the previous section we have studied the estimation $\ln \sigma$ based on i.i.d. sample. Here we will discuss the same estimation problem three special sampling schemes.  These schemes are namely (i) record values, (ii) type-II censoring, and (iii) progressive Type-II censoring. Under these sampling schemes, we derive improved estimators over the MRIE and we will observe that the results follows from the i.i.d. sampling scheme. \cite{patra2018estimating} studied the problem of estimating hazard rate of under these sampling schemes. 

\subsection{Record Values}
Various application of record model have been found in several areas such as sports analysis, hydrology, meteorology and stock market analysis. Several authors have investigated record values because of its importance. For a detail literature review in this direction we refer to \cite{ahsanullah1995introduction}, \cite{ahsanullah2001ordered} and \cite{arnold2011records}. Let $Z_1,Z_2,Z_3,\dots$ be sequence of i.i.d random variables taken form an exponential population $E(\mu,\sigma)$. For $m \ge 2$ define $u(1)=1$ and $u(m)=\min\{j|j>u(m-1),Z_j>Z_{u(m-1)}\}$, then $\{X_m=V_{u(m)}, m\ge1\}$ gives a sequence of (maximal) record statistics. 
The sequence ${u(m), m\ge 1}$ is called record times. Consider the record sample $Y_{i1},\dots,Y_{in}$ from $Exp(\theta_i,\sigma)$, $i=1,2,\dots,k.$ Then $(Y_{11},Y_{2,1},\dots,Y_{k1},S)$ be the sufficient statistics for $(\theta_1,\theta_2,\dots,\theta_k,\sigma)$, where $S=\sum_{i=1}^k(Y_{in}-Y_{i1}) \sim Gamma(k(n-1),\sigma)$ and $Y_{i1} \sim Exp(\theta_i,\sigma)$.

We have the MRIE of $\ln \sigma$ is  
$$T^R_{0}=\ln S+q_0,$$
where $q_0$ minimizes 
$$E_{\boldsymbol{\theta},1}L\left[L(\ln S+c)\right]$$
Now define $z_i=\frac{Y_{i1}}{S}, i=1,2\dots,k$ and denote $\boldsymbol{Y}=(Y_{11},Y_{21},\dots,Y_{k1})$. Then using Theorem (\ref{thmst}) we can prove that the estimator 
\begin{eqnarray}\label{rstein}
	T^R_{\zeta_0}(\boldsymbol{Y},S)=\left\{
	\begin{array}{ll}
		\ln S+\min\left\{q_0,p_0+\ln\left(\sum_{i=1}^kz_i+1\right)\right\},~~~z_i>0, i=1,\dots,k\\\\
		\ln S+q_0, ~~~~~~~~~~~~~~~~~~~~~~~~~~~~\mbox{otherwise}
	\end{array},
	\right.
\end{eqnarray}
have uniformly smaller risk than that of MRIE $T^R_{0}$, where $q_0$ and $p_0$ is given as in Theorem (\ref{thmst}).
The estimator $T^R_{\zeta_0}(\boldsymbol{Y},S)$ is non smooth. Now we will propose an estimator of $\ln \sigma$ based on
record values which dominates MRIE $T^R_{0}$. By Theorem \ref{bzthm}, the dominating estimator is obtained as 
\begin{eqnarray}
	T^R_{BZ}(\boldsymbol{Y},S)=\left\{
	\begin{array}{ll}
		\ln S+d(\boldsymbol{Z},\boldsymbol{0}),~~~~z_1>0,\dots,z_k>0\\\\
		\ln S+q_0 ,  ~~~~~~~~~\mbox{otherwise}
	\end{array}.
	\right.
\end{eqnarray}
where $d(\boldsymbol{Z},\boldsymbol{0})$ given in example (\ref{examp1}) for squared error loss function and in example (\ref{examp2}) for linex loss function. We denote $T_{01}^{*R}$ and  $T^R_{BZ1}$ be the improved estimators for squared error loss function. For linex loss function the improved estimators are denoted as $T_{02}^{*R}$ and  $T^R_{BZ2}$.

\subsubsection{Simulation study}
In this section, we compare the risk performance of the improved estimators based on record values generated from two exponential distributions. For simulation, 20,000 record samples of sizes $n=4,6,8$ are generated from two exponential distributions with location parameters $\theta_1$, $\theta_2$, and scale parameter $1$. We have presented the percentage risk improvement (PRI) with respect to the BAEE $T_0$ of the improved estimators for $n=4$ and $n=6$ and $n=8$, respectively, in Table \ref{table2}. The risk of $T_0$ is independent of $\theta_1$ and $\theta_2$ and constant. Risk values of $T_0$ are obtained as  for $n=4, 6, 8$ are  $0.183962$, $0.1065479$ and $0.07516415$ respectively. From the simulated values, we have the following observations.
\begin{itemize}
	\item [(i)] The performance of $T_{01}^{*R}$ and  $T^R_{BZ1}$ better than $T_{01}^{*}$ and  $T^R_{BZ1}$ respectively.
	\item [(ii)] The interval of improvement of $T_{01}^{*R}$ and  $T^R_{BZ1}$ larger than $T_{01}^{*}$ and  $T^R_{BZ1}$ respectively.
	\item [(iii)] The PRI of $T_{01}^{*R}$ and  $T^R_{BZ1}$ decreases slowly as $\theta_1$ and $\theta_2$ increases. 
\end{itemize}
We have a similar observation that can be made for the linex loss function. 
\begin{table}[h!]
	\centering
	\caption{Percentage risk improvement with respect to squared error loss function}
	\vspace{0.2cm}
	\label{table2}
	\begin{tabular}{|cl|cc|cc|cc|}
		\hline
		\multicolumn{2}{|c|}{$n$}                                                  & \multicolumn{2}{c|}{4}                                                       & \multicolumn{2}{c|}{6}                                                       & \multicolumn{2}{c|}{8}                                                       \\ \hline
		\multicolumn{1}{|c|}{$\theta_2$}               & \multicolumn{1}{c|}{$\theta_1$} & \multicolumn{1}{c|}{$T_{01}^{*R}$} & \multicolumn{1}{c|}{$T^R_{BZ1}$} & \multicolumn{1}{c|}{$T_{01}^{*R}$} & \multicolumn{1}{c|}{$T^R_{BZ1}$} & \multicolumn{1}{c|}{$T_{01}^{*R}$} & \multicolumn{1}{c|}{$T^R_{BZ1}$} \\ \hline
		\multicolumn{1}{|c|}{\multirow{11}{*}{0.1}} & 0.1                          & \multicolumn{1}{l|}{9.995703}          & 7.56919                             & \multicolumn{1}{l|}{6.568983}          & 5.769624                            & \multicolumn{1}{l|}{4.768522}          & 4.621218                            \\ \cline{2-8} 
		\multicolumn{1}{|c|}{}                      & 0.2                          & \multicolumn{1}{l|}{10.61596}          & 6.92789                             & \multicolumn{1}{l|}{7.048523}          & 4.693227                            & \multicolumn{1}{l|}{5.174139}          & 3.797647                            \\ \cline{2-8} 
		\multicolumn{1}{|c|}{}                      & 0.5                          & \multicolumn{1}{l|}{11.04688}          & 2.00213                             & \multicolumn{1}{l|}{7.476931}          & 2.054984                            & \multicolumn{1}{l|}{5.580389}          & 1.81209                             \\ \cline{2-8} 
		\multicolumn{1}{|c|}{}                      & 0.6                          & \multicolumn{1}{l|}{10.81651}          & 1.800788                            & \multicolumn{1}{l|}{7.327636}          & 1.37326                             & \multicolumn{1}{l|}{5.484115}          & 1.290101                            \\ \cline{2-8} 
		\multicolumn{1}{|c|}{}                      & 0.7                          & \multicolumn{1}{l|}{10.43916}          & 1.021211                            & \multicolumn{1}{l|}{7.060799}          & 0.7919498                           & \multicolumn{1}{l|}{5.288744}          & 0.880932                            \\ \cline{2-8} 
		\multicolumn{1}{|c|}{}                      & 0.8                          & \multicolumn{1}{l|}{9.962606}          & 0.3461317                           & \multicolumn{1}{l|}{6.69991}           & 0.2828253                           & \multicolumn{1}{l|}{5.012398}          & 0.420678                            \\ \cline{2-8} 
		\multicolumn{1}{|c|}{}                      & 0.9                          & \multicolumn{1}{l|}{9.397655}          & 0.2377135                           & \multicolumn{1}{l|}{6.263447}          & 0.1615921                           & \multicolumn{1}{l|}{4.673977}          & 0.062847                            \\ \cline{2-8} 
		\multicolumn{1}{|c|}{}                      & 1                            & \multicolumn{1}{l|}{8.780024}          & 0.7418085                           & \multicolumn{1}{l|}{5.778098}          & 0.5480858                           & \multicolumn{1}{l|}{4.287963}          & 0.2500554                           \\ \cline{2-8} 
		\multicolumn{1}{|c|}{}                      & 1.2                          & \multicolumn{1}{l|}{7.482334}          & 1.549307                            & \multicolumn{1}{l|}{4.742903}          & 1.171272                            & \multicolumn{1}{l|}{3.438672}          & 0.7579012                           \\ \cline{2-8} 
		\multicolumn{1}{|c|}{}                      & 1.3                          & \multicolumn{1}{l|}{6.829807}          & 1.868922                            & \multicolumn{1}{l|}{4.224752}          & 1.418521                            & \multicolumn{1}{l|}{3.004381}          & 0.9604498                           \\ \cline{2-8} 
		\multicolumn{1}{|c|}{}                      & 1.5                          & \multicolumn{1}{l|}{5.585653}          & 2.372942                            & \multicolumn{1}{l|}{3.235294}          & 1.806937                            & \multicolumn{1}{l|}{2.182459}          & 1.279523                            \\ \hline
		\multicolumn{1}{|c|}{\multirow{11}{*}{0.5}} & 0.1                          & \multicolumn{1}{l|}{11.04688}          & 2.69118                             & \multicolumn{1}{l|}{7.476931}          & 2.029748                            & \multicolumn{1}{l|}{5.580389}          & 1.809013                            \\ \cline{2-8} 
		\multicolumn{1}{|c|}{}                      & 0.2                          & \multicolumn{1}{l|}{10.81651}          & 1.187459                            & \multicolumn{1}{l|}{7.327636}          & 1.011165                            & \multicolumn{1}{l|}{5.484115}          & 1.035108                            \\ \cline{2-8} 
		\multicolumn{1}{|c|}{}                      & 0.5                          & \multicolumn{1}{l|}{9.397655}          & 2.215693                            & \multicolumn{1}{l|}{6.263447}          & 1.388422                            & \multicolumn{1}{l|}{4.673977}          & 0.8212147                           \\ \cline{2-8} 
		\multicolumn{1}{|c|}{}                      & 0.6                          & \multicolumn{1}{l|}{8.780024}          & 3.058793                            & \multicolumn{1}{l|}{5.778098}          & 2.004203                            & \multicolumn{1}{l|}{4.287963}          & 1.305989                            \\ \cline{2-8} 
		\multicolumn{1}{|c|}{}                      & 0.7                          & \multicolumn{1}{l|}{8.136173}          & 3.787176                            & \multicolumn{1}{l|}{5.266646}          & 2.543115                            & \multicolumn{1}{l|}{3.869757}          & 1.733374                            \\ \cline{2-8} 
		\multicolumn{1}{|c|}{}                      & 0.8                          & \multicolumn{1}{l|}{7.482334}          & 4.415578                            & \multicolumn{1}{l|}{4.742903}          & 3.01307                             & \multicolumn{1}{l|}{3.438672}          & 2.108569                            \\ \cline{2-8} 
		\multicolumn{1}{|c|}{}                      & 0.9                          & \multicolumn{1}{l|}{6.829807}          & 4.956729                            & \multicolumn{1}{l|}{4.224752}          & 3.42126                             & \multicolumn{1}{l|}{3.004381}          & 2.43641                             \\ \cline{2-8} 
		\multicolumn{1}{|c|}{}                      & 1                            & \multicolumn{1}{l|}{6.195272}          & 5.421653                            & \multicolumn{1}{l|}{3.719696}          & 3.774205                            & \multicolumn{1}{l|}{2.580883}          & 2.721376                            \\ \cline{2-8} 
		\multicolumn{1}{|c|}{}                      & 1.2                          & \multicolumn{1}{l|}{5.005209}          & 6.159832                            & \multicolumn{1}{l|}{2.783005}          & 4.337397                            & \multicolumn{1}{l|}{1.820591}          & 3.178925                            \\ \cline{2-8} 
		\multicolumn{1}{|c|}{}                      & 1.3                          & \multicolumn{1}{l|}{4.459724}          & 6.448651                            & \multicolumn{1}{l|}{2.369488}          & 4.557787                            & \multicolumn{1}{l|}{1.497925}          & 3.358831                            \\ \cline{2-8} 
		\multicolumn{1}{|c|}{}                      & 1.5                          & \multicolumn{1}{l|}{3.481266}          & 6.897491                            & \multicolumn{1}{l|}{1.665062}          & 4.897876                            & \multicolumn{1}{l|}{0.9633909}         & 3.637013                            \\ \hline
		\multicolumn{1}{|c|}{\multirow{11}{*}{0.7}} & 0.1                          & \multicolumn{1}{l|}{10.43916}          & 1.010742                            & \multicolumn{1}{l|}{7.060799}          & 0.785796                            & \multicolumn{1}{l|}{5.288744}          & 0.824321                            \\ \cline{2-8} 
		\multicolumn{1}{|c|}{}                      & 0.2                          & \multicolumn{1}{l|}{9.962606}          & 0.4624575                           & \multicolumn{1}{l|}{6.69991}           & 0.2081989                           & \multicolumn{1}{l|}{5.012398}          & 0.070203                            \\ \cline{2-8} 
		\multicolumn{1}{|c|}{}                      & 0.5                          & \multicolumn{1}{l|}{8.136173}          & 3.788664                            & \multicolumn{1}{l|}{5.266646}          & 2.544073                            & \multicolumn{1}{l|}{3.869757}          & 1.734094                            \\ \cline{2-8} 
		\multicolumn{1}{|c|}{}                      & 0.6                          & \multicolumn{1}{l|}{7.482334}          & 4.610197                            & \multicolumn{1}{l|}{4.742903}          & 3.141557                            & \multicolumn{1}{l|}{3.438672}          & 2.203726                            \\ \cline{2-8} 
		\multicolumn{1}{|c|}{}                      & 0.7                          & \multicolumn{1}{l|}{6.829807}          & 5.318709                            & \multicolumn{1}{l|}{4.224752}          & 3.663454                            & \multicolumn{1}{l|}{3.004381}          & 2.616951                            \\ \cline{2-8} 
		\multicolumn{1}{|c|}{}                      & 0.8                          & \multicolumn{1}{l|}{6.195272}          & 5.928763                            & \multicolumn{1}{l|}{3.719696}          & 4.117574                            & \multicolumn{1}{l|}{2.580883}          & 2.978904                            \\ \cline{2-8} 
		\multicolumn{1}{|c|}{}                      & 0.9                          & \multicolumn{1}{l|}{5.585653}          & 6.452934                            & \multicolumn{1}{l|}{3.235294}          & 4.511014                            & \multicolumn{1}{l|}{2.182459}          & 3.294357                            \\ \cline{2-8} 
		\multicolumn{1}{|c|}{}                      & 1                            & \multicolumn{1}{l|}{5.005209}          & 6.902107                            & \multicolumn{1}{l|}{2.783005}          & 4.850207                            & \multicolumn{1}{l|}{1.820591}          & 3.567731                            \\ \cline{2-8} 
		\multicolumn{1}{|c|}{}                      & 1.2                          & \multicolumn{1}{l|}{3.951458}          & 7.611995                            & \multicolumn{1}{l|}{1.997429}          & 5.388569                            & \multicolumn{1}{l|}{1.210856}          & 4.004263                            \\ \cline{2-8} 
		\multicolumn{1}{|c|}{}                      & 1.3                          & \multicolumn{1}{l|}{3.481266}          & 7.888071                            & \multicolumn{1}{l|}{1.665062}          & 5.597741                            & \multicolumn{1}{l|}{0.9633909}         & 4.174644                            \\ \cline{2-8} 
		\multicolumn{1}{|c|}{}                      & 1.5                          & \multicolumn{1}{l|}{2.657382}          & 8.313801                            & \multicolumn{1}{l|}{1.122664}          & 5.917473                            & \multicolumn{1}{l|}{0.5815305}         & 4.435515                            \\ \hline
		\multicolumn{1}{|c|}{\multirow{11}{*}{0.9}} & 0.1                          & \multicolumn{1}{l|}{9.397655}          & 0.2498145                           & \multicolumn{1}{l|}{6.263447}          & 0.1684079                           & \multicolumn{1}{c|}{4.673977}          & 0.058519                            \\ \cline{2-8} 
		\multicolumn{1}{|c|}{}                      & 0.2                          & \multicolumn{1}{l|}{8.780024}          & 1.697961                            & \multicolumn{1}{l|}{5.778098}          & 1.141567                            & \multicolumn{1}{l|}{4.287963}          & 0.6785448                           \\ \cline{2-8} 
		\multicolumn{1}{|c|}{}                      & 0.5                          & \multicolumn{1}{l|}{6.829807}          & 4.959899                            & \multicolumn{1}{l|}{4.224752}          & 3.422922                            & \multicolumn{1}{l|}{3.004381}          & 2.437712                            \\ \cline{2-8} 
		\multicolumn{1}{|c|}{}                      & 0.6                          & \multicolumn{1}{l|}{6.195272}          & 5.763112                            & \multicolumn{1}{l|}{3.719696}          & 4.004595                            & \multicolumn{1}{l|}{2.580883}          & 2.894121                            \\ \cline{2-8} 
		\multicolumn{1}{|c|}{}                      & 0.7                          & \multicolumn{1}{l|}{5.585653}          & 6.45462                             & \multicolumn{1}{l|}{3.235294}          & 4.511722                            & \multicolumn{1}{l|}{2.182459}          & 3.294943                            \\ \cline{2-8} 
		\multicolumn{1}{|c|}{}                      & 0.8                          & \multicolumn{1}{l|}{5.005209}          & 7.048864                            & \multicolumn{1}{l|}{2.783005}          & 4.952036                            & \multicolumn{1}{l|}{1.820591}          & 3.645259                            \\ \cline{2-8} 
		\multicolumn{1}{|c|}{}                      & 0.9                          & \multicolumn{1}{l|}{4.459724}          & 7.558307                            & \multicolumn{1}{l|}{2.369488}          & 5.332562                            & \multicolumn{1}{l|}{1.497925}          & 3.949792                            \\ \cline{2-8} 
		\multicolumn{1}{|c|}{}                      & 1                            & \multicolumn{1}{l|}{3.951458}          & 7.993733                            & \multicolumn{1}{l|}{1.997429}          & 5.659663                            & \multicolumn{1}{l|}{1.210856}          & 4.212918                            \\ \cline{2-8} 
		\multicolumn{1}{|c|}{}                      & 1.2                          & \multicolumn{1}{l|}{3.050316}          & 8.678725                            & \multicolumn{1}{l|}{1.374264}          & 6.17607                             & \multicolumn{1}{l|}{0.7543619}         & 4.63079                             \\ \cline{2-8} 
		\multicolumn{1}{|c|}{}                      & 1.3                          & \multicolumn{1}{l|}{2.657382}          & 8.943496                            & \multicolumn{1}{l|}{1.122664}          & 6.375266                            & \multicolumn{1}{l|}{0.5815305}         & 4.792679                            \\ \cline{2-8} 
		\multicolumn{1}{|c|}{}                      & 1.5                          & \multicolumn{1}{l|}{1.982989}          & 9.348586                            & \multicolumn{1}{l|}{0.7328501}         & 6.676806                            & \multicolumn{1}{l|}{0.3320165}         & 5.038056                            \\ \hline
	\end{tabular}
\end{table}
\subsection{Type-II censoring}
Researcher often encountered in reliability and life-testing experiments in which experimental units are either lost or removed from the experiment before failure. For example, experimental units breaks down accidentally before time in many industrial experiments; an individual withdraw from a clinical trial or the experiment may be terminated due to lack of funds. Experimenter intentionally may terminate the experiment to save time and cost associated with testing. Data obtained from such type of experiments are called censored data. One such censoring scheme is Type-II censoring. In this scheme the experimenter decides to terminate the experiment after a specified number of items $r \le n$ fail. For further details on this topic one may refer to \cite{balakrishnan2000progressive}.

Let a sample of size $n$ be drawn from an exponential distribution $E(\theta_i,\sigma)$ and the observations are available in order, that is,  $X_{i(1)} \le X_{i(2)} \le \dots \le X_{i(n)}$ for $i=1,2,\dots,k$. Here $X_{i(j)}$ is the $j^{th}$ smallest observation in a sample of $n$ observation taken from exponential $E(\theta_i,\sigma)$ population. Now consider the first $r$ ordered observations $X_{i(1)}, X_{i(2)}, \dots, X_{i(r)}$, $r\le n$, $i=1,2,\dots, k$. We consider the estimation of $\ln \sigma$ based on censored sample under bowl shaped location in variant loss function $L(t)$. Define $S=\sum_{j=1}^{k}\left[\sum_{i=1}^{r}(X_{j(i)}-X_{j(1)})+(n-r)(X_{j(n)}-X_{j(1)})\right]$. Then in this set-up $(X_{1(1)},X_{2(1)},\dots,X_{k(1)},S)$ is a minimal sufficient statistic, where for $i=1,2,\dots,k$, $X_i=nX_{i(1)}$ and $S$ follow exponential distribution $E(n\theta_i,\sigma)$ and gamma distribution $Gamma(k(n-1),\sigma)$ respectively. Consequently, the improved estimators of $\ln \sigma$ can be derived using Theorems \ref{thmst} and \ref{bzthm}.

\subsection{Progressive Type-II censoring}
Under censoring lifetimes distributions are more popular due to wide applications in science, engineering, social sciences, public health and medicine. There are several censored scheme. One important censoring scheme is  progressive Type-II censoring. Censored data are of progressively Type-II when they are censored by removing a prefixed number of surviving units when an individual unit fails.  This process continues until a fixed number of failures has occurred, at which stage the remainder of the surviving individuals are also removed/censored. For detailed one can see \cite{viveros1994interval}, \cite{balakrishnan2000progressive}, \cite{balakrishnan2007progressive}.

Now we will describe the progressive Type-II censoring scheme. The description here is similar to \cite{patra2018estimating}. Let $X_{i1},X_{2},\dots,X_{im}$ be life times of $m$  independent units placed on a life testing experiment with $X_{ij}$ following an exponential distribution $E(\theta_i,\sigma),i=1,2,\dots,k$. For $r=1,2,\dots,n$, $n \le m$, at the time of $r-$th failure, a prefixed number of $R_r$ surviving units are withdrawn from the experiment, where $R_n=m-n-R_1-R_2-\dots-R_{n-1}$. Let $X_{i1:n:m} \le X_{i2:n:m} \le \dots \le X_{in:n:m}$ be the corresponding progressive Type-II censored sample for $i=1,2,\dots,k$. We consider the estimation of $\ln \sigma$ based on progressive Type-II censored sample .
Define $S=\sum_{i=1}^k\sum_{j=1}^{n}\left[(R_j+1)(X_{ij:n:m}-X_{i1:n:m})\right]$. In this case $(X_{11:n:m},X_{21:n:m},\dots X_{k1:n:m},S)$ is a minimal sufficient statistic. Define $X_{i}=nX_{i1:n:m},i=1,2,\dots,k$. Then $X_{i}$ follows an exponential distribution $E(n\mu_{i},\sigma)$ and $S$ follows a gamma distribution $Gamma(k(n-1),\sigma)$. Consequently the improved estimators of $\ln \sigma$ can be found by Theorem \ref{thmst} and \ref{bzthm}.

\section{Conclusions} \label{sec7}
In several areas of applied statistics such as reliability engineering, molecular biology, finance, information theory, statistical physics etc., the measure of uncertainty of a probability distribution plays an important role. The Shannon's and R\'{e}nyi
entropy are the widely used measure of uncertainty. Similar to mean, standard deviation, variance and quantile, entropy is also an important
characteristic of a parametric family of distributions. In the present manuscript, we deal with the problem of estimating the entropy of several exponential distributions with respect to the bowl-shaped location invariant loss function. At first, we derived MRIE based on $S$. Now using the information contained in $(\boldsymbol{X},S)$, we have derived estimators which improve upon the MRIE of the entropy $\ln \sigma$. The techniques of \cite{stein1964}, and \cite{brewster1974improving} have been adopted to derive improved estimators. As an application, we have derived the improved estimators for squared error and linex loss functions. We have observed that the the improved estimators for (i) record values (ii) type-II censoring (iii) progressive type-II censoring can be obtaied using the results of i.i.d. sampling.  Finally we have conducted a simulation study  to compare the risk performance of the proposed estimators numerically. From the simulation, it is seen that the performance of the improved estimators is better for the record sample.  
\bibliography{censor_entropy}
\end{document}